\documentclass[12pt]{article}
\usepackage{amsfonts}

\newcommand{\sect}[1]{\setcounter{equation}{0}\section{#1}}
\newcommand{\subsect}[1]{\subsection{#1}}
\newcommand{\subsubsect}[1]{\subsubsection{#1}}
\renewcommand{\theequation}{\arabic{section}.\arabic{equation}}

\def\be{\begin{equation}}
\def\ee{\end{equation}}
\def\bea{\begin{eqnarray}}
\def\eea{\end{eqnarray}}

\def\1{\'{\i}}
\def\R{{\mathbb R}}

\def\buno{c_1}
\def\bdos{c_3}
\def\btres{c_2}
\def\bseis{c_4}
\def\bsiete{c_6}
\def\bdiez{c_{5}}

\def\aone{a_+}
\def\atwo{a_-}
\def\athree{b_+}
\def\afour{b_-}
\def\afive{a}
\def\asix{b}

\def\fama{${\mbox{I}_+}$}
\def\famb{${\mbox{I}_-}$}
\def\famc{${\mbox{II}}$}

\def\K{J_3}
\def\P{J_+}
\def\T{J_-}
\def\M{I}

\def\mmm{M}

\def\ssl{{gl}(2)}

\begin{document}

\thispagestyle{empty}

 \
\hfill math.QA/9806149

\hfill To appear in {\em J.Phys.A}

\
\vspace{1cm}

\begin{center}
{\LARGE{\bf{Multiparametric quantum $gl(2)$:}}}

{\LARGE{\bf{Lie bialgebras, quantum $R$-matrices and}}}

{\LARGE{\bf{non-relativistic limits}}}
\end{center}

\bigskip\bigskip

\begin{center} Angel Ballesteros$^\dagger$, Francisco J.
Herranz$^\dagger$  and Preeti Parashar$^{\dagger\ddagger}$
\end{center}

\begin{center} {\it {$^\dagger$Departamento de F\1sica, Universidad
de Burgos} \\   Pza. Misael Ba\~nuelos,
E-09001 Burgos, Spain}
\end{center}

\begin{center} {\it {$^\ddagger$Departamento de F\1sica Te\'orica,
Universidad Aut\'onoma de Madrid} \\   Cantoblanco,
E-28049 Madrid, Spain}
\end{center}

\bigskip\bigskip

\begin{abstract}
Multiparametric quantum deformations of $gl(2)$ are studied through a
complete classification of $gl(2)$ Lie bialgebra structures. From them,
the non-relativistic limit leading to harmonic oscillator Lie
bialgebras is implemented by means of a contraction procedure. New
quantum deformations of $gl(2)$ together with their associated quantum
$R$-matrices are obtained and other known quantizations are recovered
and classified. Several connections with integrable models are outlined.
\end{abstract}

\bigskip\bigskip
\bigskip\bigskip



\newpage


\sect{Introduction}

A two-parametric quantum deformation of
$gl(2)$ has been proven in \cite{MRplb}
to provide the quantum group symmetry of the spin
$1/2$ XXZ Heisenberg chain with twisted periodic boundary
conditions \cite{AGR,PS}. In this context, the central generator  
$\M$ of the $gl(2)$ algebra plays an essential role in the algebraic 
introduction of the twisted boundary terms of the spin Hamiltonian
through a deformation induced from the exponential of the classical
$r$-matrix $r=J_3\wedge I$. This seems not to be an isolated example,
since the general construction introduced in
\cite{FLR} establishes a correspondence between models with twisted
boundary conditions (see references therein) and
multiparametric Reshetikhin twists
\cite{Res} in which the Cartan subalgebra is enlarged with a
(cohomologically trivial) central generator.

From another different physical point of view, $gl(2)$ can be
also considered as the natural relativistic analogue of the
one-dimensional harmonic oscillator algebra \cite{Al}. The latter (which
is a non-trivial central extension of the (1+1) Poincar\'e algebra) can
be obtained from
$gl(2)$ (which is a trivial central extension of $sl(2,\R)\equiv so(2,1)$)
through a generalized In\"on\"u--Wigner contraction, that can be
interpreted as the algebraic transcription of the non-relativistic limit
connecting both kinematics. The direct applicability of quantum algebras
in the construction of completely integrable many-body systems through
the formalism given in
\cite{orl} (that precludes the use of any transfer matrix
technique by making use directly of the Hopf algebra axioms) suggests that
a systematic study of quantum
$gl(2)$ algebras would be related to the definition of integrable
systems consisting in long-range interacting relativistic oscillators (see
\cite{chain} for the construction of non-relativistic oscillator
chains). Finally, note also that a
$gl(2)$ induced deformation of the Schr\"odinger algebra has been
recently used to construct a discretized version of the (1+1)
Schr\"odinger equation on a uniform time lattice \cite{Sch}.

So far, much attention has been
paid to quantum
$GL(2)$ groups and their
classifications \cite{standb}--\cite{Kuper} but a
fully general and explicit description of quantum $gl(2)$ algebras is
still lacking, although partial results can be already found in the
literature
\cite{wess}--\cite{Preeti}.
Such a systematic approach to quantum $gl(2)$ algebras is the aim of the
present paper, and the underlying Lie bialgebra structures and
classical $r$-matrices will be shown to contain all the essential
information characterizing different quantizations. In section 2, $gl(2)$
Lie bialgebras are fully obtained and classified into two multiparametric
and inequivalent families. Their contraction to the harmonic
oscillator Lie bialgebras is performed in section 3 by
introducing a multiparameter generalization of the Lie bialgebra
contraction theory \cite{LBC} that allows us to perform the
non-relativistic limit. Among the quantum deformations of the harmonic
oscillator algebra whose Lie bialgebras are obtained, we find the one
introduced in
\cite{GS} in the context of link invariants.
Finally, an extensive study of the
quantizations of $gl(2)$ Lie bialgebras is given in section 4.
New quantum algebras, deformed Casimir operators and quantum
$R$-matrices are obtained and known
results scattered through the literature are easily derived from the
classification presented here. In particular, the quantum  algebra
corresponding to the quantum group $GL_{h,q}(2)$ \cite{Kuper} is
constructed (recall that this is a natural superposition of both standard
and non-standard deformations). Quantum symmetry algebras of the twisted
XXZ and XXX models are identified, and it is shown how new quantum
$gl(2)$ invariant spin chains can be systematically obtained from the
new multiparametric deformations that have been introduced.


\sect{$gl(2)$ Lie bialgebras}

A Lie bialgebra $(g,\delta)$ is a Lie
algebra $g$ endowed with a map $\delta:g\to g\otimes g$ (the
cocommutator) that fulfils two conditions:

\noindent i) $\delta$ is a 1-cocycle, i.e.,
\be
\delta([X,Y])=[\delta(X),\, 1\otimes Y+ Y\otimes 1] +
[1\otimes X+ X\otimes 1,\, \delta(Y)]  \qquad \forall X,Y\in
g.
\label{ba}
\ee
\noindent ii) The dual map $\delta^\ast:g^\ast\otimes g^\ast \to
g^\ast$ is a Lie bracket on $g^\ast$.

A Lie bialgebra $(g,\delta)$ is called a coboundary Lie bialgebra  if
there exists an element $r\in g\wedge g$ (the classical $r$-matrix),
such that
\be
\delta(X)=[1\otimes X + X \otimes 1,\,  r]  \qquad  \forall X\in
g.
\label{bb}
\ee
When the $r$-matrix is a skewsymmetric solution of the classical
Yang--Baxter equation (YBE)  we shall
say that $(g,\delta (r))$ is a {\em non-standard} (or triangular) Lie
bialgebra, while when it is a skewsymmetric solution of the modified
classical  YBE  we shall have a {\em standard}
one. On the other hand, two Lie bialgebras $(g,\delta)$  and
$(g,\delta')$ are said to be equivalent if there exists an automorphism
$O$ of $g$ such that
$\delta'=(O\otimes O)\circ\delta\circ O^{-1}$.

Let us now consider the $gl(2)$ Lie algebra
\be
[\K,\P]=2\P \qquad [\K,\T]=-2\T \qquad [\P,\T]=\K \qquad [\M,\cdot\,]=0 .
\label{bc}
\ee
Notice that  $\ssl=sl(2,\R)\oplus u(1)$
where $\M$ is the central generator.
The second order Casimir is
\be
{\cal C}=\K^2 + 2 \P \T + 2 \T \P.
\label{bd}
\ee
The most general cocommutator  $\delta:\ssl \rightarrow \ssl
\otimes\ssl$  will be a linear combination (with real coefficients)
\be
\delta(X_i)=f_i^{jk}\,X_j\wedge X_k
\label{be}
\ee
of skewsymmetric products of  the generators $X_l$ of $\ssl$. Such a
completely general cocommutator has to be computed by firstly  imposing
the cocycle condition.
This leads to the following six-parameter
$\{\aone,\atwo,\athree,\afour,\afive,\asix\}$ (pre)cocommutator:
\bea
&&\delta(\K)=\aone \K\wedge \P+ \atwo \K\wedge \T + \athree \P\wedge \M
+ \afour \T\wedge \M\cr
&&\delta(\P)=\afive \K\wedge \P -\frac {\afour}2 \K\wedge \M +
\atwo \P\wedge \T + \asix \P\wedge \M\cr
&&\delta(\T)=\afive \K\wedge \T -\frac {\athree}2 \K\wedge \M
- \aone \P\wedge \T - \asix \T\wedge \M\cr
&&\delta(\M)=0.
\label{bf}
\eea

Afterwards, Jacobi identities have to be imposed onto
$\delta^\ast:\ssl^\ast \otimes \ssl^\ast \rightarrow \ssl^\ast $
in order to guarantee that a  Lie bracket is defined through this map.
Thus we obtain the following set of equations:
\be
 \aone\asix - \athree\afive=0\qquad
 \aone\afour+ \atwo\athree=0\qquad
 \atwo\asix + \afour\afive=0 .
\label{bg}
\ee

The next step is to find out the  Lie  bialgebras
defined by (\ref{bf}) and (\ref{bg}) that come from classical
$r$-matrices. Let us consider an arbitrary skewsymmetric element of
$\ssl \wedge \ssl$:
\be
r=\buno \K\wedge \P  + \btres  \K\wedge \T + \bdos  \K\wedge \M
+ \bseis \P\wedge \M + \bdiez \T\wedge \M + \bsiete  \P\wedge \T .
\label{da}
\ee
The corresponding   Schouten bracket reads:
\bea
&&[[r,r]]= ({\bsiete}^2 - 4\buno\btres) \K\wedge \P\wedge \T +
(\bseis\bsiete - 2
\buno\bdos) \K\wedge \P\wedge \M  \cr
&&\qquad\qquad + (2 \bdos\btres + \bsiete\bdiez) \K\wedge \T\wedge \M +
2 (\btres\bseis +  \buno\bdiez) \P\wedge \T\wedge \M ,
\label{db}
\eea
and the   modified classical YBE  will be satisfied
provided
\be
\bseis\bsiete - 2
\buno\bdos = 0\qquad 2 \bdos\btres + \bsiete\bdiez  = 0
\qquad \btres\bseis +  \buno\bdiez = 0 .
\label{db3}
\ee
These equations  map exactly
onto the conditions (\ref{bg}) (obtained from the Jacobi identities)
under the following identification of the parameters:
\bea
&&\aone = 2 \buno\qquad \atwo = - 2 \btres\qquad \athree = 2 \bseis \cr
&&\afour = -2 \bdiez\qquad \afive = - \bsiete\qquad \asix = - 2 \bdos .
\label{dc1}
\eea
Therefore all  Lie bialgebras
associated with   $\ssl$ are coboundaries and the most general
$r$-matrix (\ref{da}) can be written in terms of the $a$'s and $b$'s
parameters:
\be r=\frac 12(\aone \K\wedge \P - \atwo  \K\wedge \T  - \asix  \K\wedge
\M + \athree \P\wedge \M - \afour \T\wedge \M - 2\afive  \P\wedge \T ).
\label{dda}
\ee
Under these conditions, the Schouten bracket  reduces to
\be
[[r,r]]= (\bsiete^2 - 4 \buno\btres) \K\wedge \P\wedge \T =
(\afive^2+\aone\atwo) \K\wedge \P\wedge \T ,
\label{ddb}
\ee
so that  it allows us to distinguish between standard
($\afive^2+\aone\atwo\ne 0)$  and non-standard
($\afive^2+\aone\atwo=0)$  Lie bialgebras.

On the other hand, the only element  $\eta\in \ssl
\otimes \ssl$  being
$Ad^{\otimes 2}$-invariant is given by
\be
\eta = \tau_1   (\K\otimes \K + 2 \T\otimes \P + 2 \P\otimes \T)
+\tau_2  \M\otimes \M ,
\label{dd}
\ee
where $\tau_1$ and $\tau_2$ are arbitrary parameters. Since $r'=r+ \eta$
will generate the same bialgebra as $r$,  the element
$\eta$ will  relate   non-skewsymmetric $r$-matrices with
skewsymmetric ones.

Let us now explicitly solve the equations (\ref{bg});
we find three disjoint families:

\noindent
$\bullet$ {\it Family \fama}:

 \noindent
Standard: $\{\aone\neq 0,\
\atwo,\ \athree,\ \afour=- {\atwo\athree}/{\aone},\ \afive,\
\asix= {\athree\afive}/{\aone}\}$ and
$\afive^2+\aone\atwo\ne 0$.

 \noindent
Non-standard: $\{\aone\neq 0,\
\atwo=- \afive^2/\aone,\ \athree,\ \afour={\athree\afive^2}/{\aone^2},\
\afive,\
\asix= {\athree\afive}/{\aone}\}$.

 \noindent
$\bullet$ {\it Family \famb}:

 \noindent
Standard: $\{\aone= 0,\
\atwo\ne 0,\ \athree=0,\ \afour,\ \afive \ne 0,\
\asix=- {\afour\afive}/{\atwo} \}$.

 \noindent
Non-standard: $\{\aone= 0,\
\atwo\ne 0,\ \athree=0,\ \afour,\ \afive=0,\
\asix=0 \}$.

 \noindent
$\bullet$ {\it Family \famc}:

 \noindent
Standard:  $\{\aone= 0,\
\atwo= 0,\ \athree=0,\ \afour=0,\ \afive\ne 0,\  \asix \}$.

 \noindent
Non-standard: $\{\aone= 0,\
\atwo= 0,\ \athree ,\  \afour,\ \afive=0,\
\asix  \}$.

This classification can be simplified by taking into account the
following  automorphism  of  $\ssl$:
\be
\P \rightarrow \T \qquad  \T \rightarrow \P  \qquad \K \rightarrow - \K
\qquad \M \rightarrow \M
\label{de}
\ee
which leaves the   Lie brackets (\ref{bc}) invariant. This map can be
implemented at a Lie bialgebra level onto
(\ref{bf}), and it implies a transformation of the deformation
parameters of the form
\be
 \aone \rightarrow \atwo\quad
 \atwo \rightarrow \aone\quad
 \athree \rightarrow -\afour\quad
 \afour \rightarrow -\athree\quad
 \afive \rightarrow - \afive\quad \asix \rightarrow - \asix .
\label{df}
\ee
The Jacobi identities (\ref{bg}) and the
classical $r$-matrix (\ref{dda}) are invariant under the
automorphism defined by (\ref{de}) and (\ref{df}).
Therefore, the family {\famb} is
included within {\fama} provided $\atwo=0$. Hence we shall
consider only the two families {\fama} and {\famc}, whose explicit
cocommutators and $r$-matrices are written in Table 1. Note that the
central generator $\M$ always has a vanishing cocommutator.

\newpage

\bigskip

{\footnotesize

 \noindent
{{\bf Table 1.} $gl(2)$ Lie bialgebras.}
\smallskip

\noindent
\begin{tabular}{lll}
\hline
 &\multicolumn{2}{c}{Family \fama}\\[0.2cm]
&Standard&Non-standard\\[0.1cm]
&$(\aone\ne 0,\atwo,\athree,\afive$ and $\afive^2+\aone \atwo\ne0$)
&$(\aone\ne 0,\athree,\afive)$\\[0.2cm]
\hline
$r$& $\frac 12 (\aone \K\wedge \P -\atwo  \K\wedge\T -
\frac{\athree\afive}{\aone}
\K\wedge \M$& $\frac 12 (\aone \K\wedge \P + \frac{\afive^2}{\aone}
\K\wedge\T -
\frac{\athree\afive}{\aone}
\K\wedge \M$\\[0.1cm]
 &$\quad +\athree \P\wedge \M + \frac{\atwo\athree}{\aone} \T\wedge
\M  - 2 \afive \P\wedge\T)$&   $\quad +\athree \P\wedge \M -
\frac{\athree\afive^2}{\aone^2} \T\wedge \M  - 2 \afive \P\wedge\T)$
\\[0.3cm]
$\delta(\K)$&$-(\aone \P + \atwo \T)\wedge  \K
+\athree (\P -\frac{\atwo}{\aone}\T) \wedge \M$&
$-\aone (\P -\frac{\afive^2}{\aone^2}  \T)\wedge \K
+\athree (\P +\frac{\afive^2}{\aone^2}\T) \wedge \M$ \\[0.1cm]
$\delta(\P)$&$(\afive \K - \atwo \T)\wedge \P
+ \frac{\athree}{ \aone}( \afive \P+\frac{\atwo}{2}\K)\wedge \M$
&$ \afive (  \K + \frac{\afive}{\aone} \T)\wedge \P
+ \frac{\athree\afive}{\aone} (\P - \frac{\afive}{2\aone}\K)\wedge
\M
$\\[0.1cm]
 $\delta(\T)$&$ (\afive \K  - \aone \P)\wedge \T
-\frac{\athree}2(\K + \frac{2\afive}{\aone}\T)\wedge \M$
&$ (\afive \K  - \aone \P)\wedge \T
 -\frac{\athree}2(\K +  \frac{2\afive}{\aone} \T)\wedge \M$\\[0.1cm]
$\delta(\M)$&$0$&$0$\\[0.2cm]
\hline
 &\multicolumn{2}{c}{Family \famc}\\[0.2cm]
&Standard&Non-standard\\[0.1cm]
&$(\afive\ne 0,\asix)$
&$(\athree,\afour,\asix)$\\[0.2cm]
\hline
$r$& $   - \frac 12  \asix \K\wedge \M - \afive \P\wedge\T $&
 $-\frac 12 (\asix \K - \athree \P + \afour \T)\wedge \M $\\[0.3cm]
$\delta(\K)$&$0$&$(\athree \P + \afour \T)\wedge \M$\\[0.1cm]
$\delta(\P)$&$-\afive \P\wedge \K + \asix \P\wedge \M$
&$- (\frac 12\afour \K - \asix \P)\wedge \M$\\[0.1cm]
$\delta(\T)$&$-\afive \T\wedge \K - \asix \T\wedge \M$
&$-(\frac 12 \athree \K + \asix \T)\wedge \M$\\[0.1cm]
$\delta(\M)$&$0$&$0$\\[0.2cm]
 \hline
\end{tabular}}

\bigskip


\subsect{$GL(2)$ Poisson--Lie groups}

It is well-known \cite{Drlb} that  when a Lie bialgebra $(g,\delta)$ is a
coboundary  one with classical $r$-matrix $r=\sum_{i,j}
r^{ij}X_i\otimes X_j$, the Poisson--Lie
bivector $\Lambda$ linked to it is given by the so called Sklyanin bracket
\be
\Lambda=\sum_{i,j} r^{ij}(X_i^L\otimes X_j^L -
X_i^R\otimes X_j^R )
\label{Sklyanin}
\ee
where $X_i^L$ and $X_j^R$  are  left and right invariant vector
fields on the Lie group $G=\mbox{Lie}(g)$.
We have just found that  all  $gl(2)$ Lie bialgebras are coboundary ones;
therefore, we can deduce their corresponding  Poisson--Lie groups  by
means of the Sklyanin bracket (\ref{Sklyanin}) as follows.

The $2\times2$ fundamental
representation
$D$ of the
$gl(2)$ algebra (\ref{bc}) is:
\bea
&&D(\K)=\left(\begin{array}{cc}
1& 0 \\ 0&-1
\end{array}\right)\qquad 
D(\M)=\left(\begin{array}{cc}
1&0\\ 0&1
\end{array}\right) \cr
&&D(\P)=\left(\begin{array}{cc}
0&1 \\ 0&0
\end{array}\right) \qquad
D(\T)=\left(\begin{array}{cc}
0& 0 \\ 1&0
\end{array}\right) .
\label{frep}
\eea
By using this representation, a group element of
$GL(2)$ can be written as
\be
T = e^{\theta_- D(\T)}e^{\theta
D(\M)}e^{\theta_3 D(\K)}
e^{\theta_+ D(\P)} =\left(\begin{array}{cc}
e^{\theta + \theta_3}&e^{\theta + \theta_3}\theta_+  \\e^{\theta +
\theta_3}\theta_- &e^{\theta + \theta_3}\theta_-\theta_+ + e^{\theta -
\theta_3} \end{array}\right).
\label{tmatrix}
\ee
Now, left and right invariant $GL(2)$ vector fields can be obtained:
\bea
&&X_{\K}^L=\partial_{\theta_3}-2\theta_+\partial_{\theta_+}
\qquad\quad X_{\M}^L=\partial_{\theta}\cr
&&X_{\P}^L=\partial_{\theta_+}\qquad\quad
X_{\T}^L=
\theta_+\partial_{\theta_3}-\theta_+^2\partial_{\theta_+}+e^{-2\theta_3}
\partial_{\theta_-}
\label{izfield}
\eea
\bea
&&X_{\K}^R=\partial_{\theta_3}-2\theta_-\partial_{\theta_-}
\qquad\quad X_{\M}^R=\partial_{\theta}\cr
&&X_{\P}^R=\theta_-\partial_{\theta_3}-\theta_-^2\partial_{\theta_-}
+e^{-2\theta_3}
\partial_{\theta_+}
 \qquad\quad
X_{\T}^R=\partial_{\theta_-} .
\label{defield}
\eea
By substituting (\ref{izfield}), (\ref{defield}) and the  classical
$r$-matrix (\ref{dda}) within  the Sklyanin bracket (\ref{Sklyanin}) we
obtain the following Poisson--Lie brackets among the (local) coordinates
$\{\theta_-,\theta_+,\theta,\theta_3\}$:
\bea
&&\{\theta_+,\theta_3\} =- \afive \theta_+ + \frac{\atwo}{2}
\theta_+^2 -
\frac{\aone}{2} (1 - e^{-2 \theta_3})  \cr
&&\{\theta_-, \theta_3\} = -\afive \theta_- + \frac{\aone}{2}
\theta_-^2 -
\frac{\atwo}{2} (1 - e^{-2 \theta_3})  \cr
&&\{\theta_+, \theta_-\} = (\atwo \theta_+ -\aone \theta_- ) e^{-2
\theta_3} \cr
&&\{\theta_+, \theta\} = \asix \theta_+ + \frac{\afour}{2}\theta_+^2 +
\frac{\athree}{2} (1 - e^{-2 \theta_3})   \label{buena} \\
&&\{\theta_-, \theta\} = -\asix \theta_- + \frac{\athree}{2} \theta_-^2 +
\frac{\afour}{2} (1 - e^{-2 \theta_3})  \cr
&&\{\theta_3, \theta\} = - \frac{1}{2} (\athree\theta_- +\afour\theta_+).
\nonumber
\eea
By imposing Jacobi identities onto (\ref{buena}), conditions (\ref{bg})
restricting the space of Lie bialgebras are recovered.  On the other
hand, from  (\ref{buena}) and (\ref{bg}) it is immediate to obtain
explicitly the Poisson--Lie groups associated to the families of $gl(2)$
Lie bialgebras written in Table 1.

We recall that a  classification of Poisson--Lie
structures on the group $GL(2)$ was carried out by Kupershmidt  in
\cite{Kuper}, where quantum group structures on $GL(2)$ were also
analysed. The relationship between (\ref{buena}) and such a
classification can be explored by  writing the matrix
$T$ (\ref{tmatrix}) as
\be
T=\left(\begin{array}{cc}
A & B \\ C & D \end{array}\right) .
\label{mt}
\ee
Starting from the Poisson brackets (\ref{buena}), the  quadratic
Poisson brackets between $\{A,B,C,D\}$ are obtained:
\bea
&&\{A,C\}=(a+b) AC-\left(\frac{a_+ +b_+}{2}\right) C^2 +\left(\frac{a_- -
b_-}{2}\right)(A^2+BC-AD)\cr 
&&\{A,B\}=(a-b) AB-\left(\frac{a_- +b_-}{2}\right) B^2
+\left(\frac{a_+ - b_+}{2}\right)(A^2+BC-AD)\cr 
&&\{B,D\}=(a+b) BD-\left(\frac{a_+ +b_+}{2}\right) (D^2+BC-AD)
+\left(\frac{a_- - b_-}{2}\right) B^2 \cr
&&\{C,D\}=(a-b) CD+\left(\frac{a_+ -b_+}{2}\right) C^2-\left(\frac{a_- +
b_-}{2}\right) (D^2+BC-AD) \cr
&&\{A,D\}=2a BC -\left(\frac{a_+ +b_+}{2}\right) CD+\left(\frac{a_+ - 
b_+}{2}\right) AC - \left(\frac{a_- + b_-}{2}\right) BD \cr
&&\qquad\qquad  + \left(\frac{a_- -
b_-}{2}\right) AB \cr 
&&\{B,C\}=2b BC -\left(\frac{a_+ +b_+}{2}\right) CD-\left(\frac{a_+ - 
b_+}{2}\right) AC + \left(\frac{a_- + b_-}{2}\right) BD   \cr
&&\qquad\qquad +\left(\frac{a_- -
b_-}{2}\right) AB .\label{kkpp}
\eea
Therefore the Poisson structures  given in \cite{Kuper} can be
completely embedded within (\ref{kkpp}) provided the following 
identification is imposed:
\bea
&&r = \afive + \asix  \qquad s = -(\aone + \athree)/2 
\qquad v = \asix - \afive \cr
&&u = \athree - \aone  \qquad w = 2 \atwo = 2 \afour.
\label{ident}
\eea
Here, $r,s,u,v$ and $w$ are the parameters arising in Kupershmidt's
classification and $\{A,B,C,D\}$ are the corresponding generators 
(note that we have
used capitals for the latter in order to avoid confusion with  the
$gl(2)$ Lie bialgebra parameters). From (\ref{ident}) we conclude that
Lie bialgebras having
$\atwo
\neq \afour$ have no counterpart in \cite{Kuper}. For
the remaining cases, (\ref{ident}) gives a straightforward 
correspondence between the quantum algebras that will be obtained after
quantization and the quantum
$GL(2)$ groups described in \cite{Kuper}.


\sect{Harmonic oscillator Lie bialgebras through contractions}

The   $gl(2)$ algebra is isomorphic to the relativistic
oscillator algebra introduced in \cite{Al} and its natural
non-relativistic limit is the harmonic oscillator algebra $h_4$. Both
algebras are related by means of a generalized In\"on\"u--Wigner
contraction
\cite{Weimar}. If we define
\be
A_+ = \varepsilon \P\qquad A_- = \varepsilon \T\qquad
N = \frac 12 (\K + \M)\qquad  \mmm = \varepsilon^2 \M ,
\label{na}
\ee
the limit $\varepsilon \to 0$ of the Lie
brackets obtained from (\ref{bc}) yields the oscillator algebra $h_4$
\be
[N,A_+]=A_+ \qquad [N,A_-]=-A_- \qquad
[A_-,A_+]=\mmm
\qquad [\mmm,\cdot\,]=0,
\label{nb}
\ee
and the parameter can be interpreted as $\varepsilon=1/c$, with $c$ being
the speed of light.

In the sequel we work out the contractions from the multiparameter
$gl(2)$ bialgebras written in Table 1 to multiparameter $h_4$
bialgebras. The Lie bialgebra contraction (LBC) approach was
introduced in  \cite{LBC} for a single deformation parameter.
In order to perform an LBC we
need two maps: the Lie
algebra transformation (an In\"on\"u--Wigner
contraction as (\ref{na})) together with a mapping on the
initial deformation parameter $a$
\be
a=\varepsilon^{n} a'
\label{nbb}
\ee
 where $n$ is any real number and $a'$ is the
contracted deformation parameter. The convergency  of the classical
$r$-matrix and the cocommutator $\delta$ under the limit
$\varepsilon \to 0$ have to be analyzed
separately, since  starting 
from a coboundary bialgebra, the LBC can lead  to another coboundary
bialgebra (both $r$ and $\delta$ converge) or can   produce  a
non-coboundary bialgebra  ($r$ diverges but $\delta$
converges). In other words, we have to find out the
minimal value of the number $n$  such that $r$ converges, 
the minimal value of $n$ such that $\delta$ converges,
and finally to compare both of them \cite{LBC}.

In the sequel, we show that the LBC method can be applied to 
multiparameter Lie bialgebras by considering a different map (\ref{nbb})
for each deformation parameter. Let us describe this procedure by
contracting, for instance,  the non-standard family II given in Table 1.

 First, we
analyze the classical $r$-matrix. We consider the following maps
\be
\athree= 2\varepsilon^{n_+}\beta_+ \qquad
\afour=- 2\varepsilon^{n_-}\beta_-\qquad
\asix=-\varepsilon^{n}\vartheta  
\label{nbc}
\ee
where $\beta_+$, $\beta_-$, $\vartheta$ are the contracted
deformation parameters, and $n_+$, $n_-$, $n$ are real numbers to be
determined by imposing the convergency of   $r$ under the non-relativistic
limit. We introduce the Lie algebra contraction (\ref{na}) and the maps
(\ref{nbc})  in the non-standard classical $r$-matrix of the family
II:
\bea
&&r=-\frac 12 (\asix \K - \athree \P + \afour \T)\wedge \M\cr
&&\quad 
=\frac 12\left(\varepsilon^{n}\vartheta (2 N - \mmm
\varepsilon^{-2}) +  2\varepsilon^{n_+}\beta_+ A_+ \varepsilon^{-1}
+  2\varepsilon^{n_-}\beta_- A_- \varepsilon^{-1}\right)
\wedge \mmm \varepsilon^{-2}\cr
&&\quad =\left(\varepsilon^{n-2}\vartheta  N
+ \varepsilon^{n_+ -3}\beta_+ A_+  + 
 \varepsilon^{n_- -3}\beta_- A_-  \right)
\wedge \mmm .
\label{nbd}
\eea
Thus the minimal values of $n$, $n_+$, $n_-$ which allow $r$ to
converge under the limit $\varepsilon \to 0$ are given by
\be
\qquad n= 2\qquad n_+= 3\qquad  n_-= 3 ,
\label{nbe}
\ee 
and the contracted $r$-matrix turns out to be
\be
r= ( \vartheta  N+  \beta_+ A_+  + \beta_- A_-   )\wedge \mmm .
\label{nbf}
\ee

Likewise we analyze the convergency of $\delta$:
\bea
&&\delta(N)=\frac 12 \left(\delta(\K) + \delta(\M)\right)=
\frac 12  (\athree \P + \afour \T)\wedge \M\cr
&&\qquad =\frac 12 \left(2\varepsilon^{n_+}\beta_+ A_+
\varepsilon^{-1} - 2\varepsilon^{n_-}\beta_- A_-
\varepsilon^{-1}\right)\wedge \mmm \varepsilon^{-2}\cr
&&\qquad =\left( \varepsilon^{n_+-3}\beta_+ A_+
  -  \varepsilon^{n_- -3}\beta_- A_- \right)\wedge \mmm ,
\label{nbh}
\eea
\bea
&&\delta(A_+)=\varepsilon\delta(\P)=
- \varepsilon (\frac 12\afour \K - \asix \P)\wedge \M\cr
&&\qquad = \varepsilon \left(\varepsilon^{n_-}\beta_- (2 N - \mmm
\varepsilon^{-2})  -\varepsilon^{n}\vartheta
A_+\varepsilon^{-1}\right) \wedge
\mmm \varepsilon^{-2}\cr
&&\qquad  =\left(2\varepsilon^{n_- -1}\beta_-   N  
-\varepsilon^{n-2}\vartheta A_+ \right) \wedge
\mmm ,
\label{nbi}
\eea
\bea
&&\delta(A_-)=\varepsilon\delta(\T)=
-\varepsilon(\frac 12 \athree \K + \asix \T)\wedge \M\cr
&&\qquad = -\varepsilon \left(\varepsilon^{n_+}\beta_+ (2 N - \mmm
\varepsilon^{-2})  -\varepsilon^{n}\vartheta
A_-\varepsilon^{-1}\right) \wedge
\mmm \varepsilon^{-2}\cr
&&\qquad  =-\left(2\varepsilon^{n_+ -1}\beta_+   N  
-\varepsilon^{n-2}\vartheta A_- \right) \wedge
\mmm ,
\label{nbj}
\eea
and, obviously, $\delta(\mmm)=0$. Hence the minimal values
of $n$, $n_+$, $n_-$ which ensure the convergency of $\delta$ under
the limit $\varepsilon\to 0$ read
\be
\qquad n= 2\qquad n_+= 3\qquad 
n_-= 3 ,
\label{nbk}
\ee 
and the contracted cocommutator reduces to
\bea
 &&\delta(N)= (  \beta_+ A_+
  -   \beta_- A_-  )\wedge \mmm \qquad
\delta(\mmm)=0\cr
&&\delta(A_+)=    - \vartheta A_+   \wedge \mmm
\qquad
\delta(A_-) =   
 \vartheta A_-   \wedge \mmm .
\label{nbl}
\eea
Therefore,   in this case, the resulting contracted
bialgebra is a coboundary one, as the contraction exponents
coming from (\ref{nbe}) and (\ref{nbk}) coincide.

The remaining $gl(2)$ families of bialgebras can be contracted by
following the LBC approach, and  all the resulting contracted
bialgebras are coboundaries. Transformations of the deformation
parameters  for the  LBCs  of the families
{\fama} and {\famc} read

\begin{tabular}{lllll}
 {\fama} Standard: &
$ \aone=\varepsilon \alpha_+$&$
 \atwo=-\varepsilon^3\beta_-$&$ \athree=-\varepsilon\alpha_+$&$
\afive= \varepsilon^2\vartheta$\cr
 {\fama} Non-standard: &
 $\aone=\varepsilon \alpha_+$&$ \athree=-\varepsilon\alpha_+$&$
\afive= \varepsilon^2\vartheta$&\cr
 {\famc} Standard: &
$\afive= -\varepsilon^2\xi$&$ \asix=-\varepsilon^2\vartheta$& &\cr
 {\famc} Non-standard: &
$\athree= 2\varepsilon^3\beta_+$&$
\afour=- 2\varepsilon^3\beta_-$&$
\asix=-\varepsilon^2\vartheta$ &\cr
\end{tabular}

If we apply these maps together with  (\ref{na})  to the
$\ssl$ Lie bialgebras displayed in Table 1 and we take the limit
$\varepsilon\to 0$, then the oscillator Lie bialgebras given in Table 2
are derived. We stress that the LBC procedure just described can be
applied in a similar way to any arbitrary multiparametric Lie
bialgebra.

\newpage

\bigskip

{\footnotesize

 \noindent
{{\bf Table 2.} Harmonic oscillator $h_4$ bialgebras via contraction
from $\ssl$.}
\smallskip

\noindent
\begin{tabular}{lll}
\hline
 &\multicolumn{2}{c}{Family \fama}\\[0.2cm]
&Standard&Non-standard\\[0.1cm]
&$(\alpha_+\ne 0,\vartheta,\beta_-$ and $\vartheta^2-
\alpha_+\beta_-\ne0$) &$(\alpha_+\ne 0,\vartheta)$\\[0.2cm]
\hline
$r$& $\alpha_+ N\wedge A_+ +\vartheta (N\wedge\mmm - A_+\wedge A_-)$&
$\alpha_+ N\wedge A_+ +\vartheta (N\wedge\mmm - A_+\wedge A_-)$\\[0.1cm]
& $\quad +\beta_- A_-\wedge \mmm$&$\quad +({\vartheta^2}/{\alpha_+})
A_-\wedge
\mmm$\\[0.3cm]
$\delta(N)$&$\alpha_+ N\wedge A_+ - \beta_- A_-\wedge \mmm$&
$\alpha_+ N\wedge A_+ - ({\vartheta^2}/{\alpha_+}) A_-\wedge \mmm$
\\[0.1cm]
$\delta(A_+)$&$0$
&$0$\\[0.1cm]
 $\delta(A_-)$&$\alpha_+(N\wedge \mmm - A_+\wedge A_-) +2 \vartheta
A_-\wedge \mmm$ &$\alpha_+(N\wedge \mmm - A_+\wedge A_-) +2 \vartheta
A_-\wedge \mmm$\\[0.1cm]
$\delta(\mmm)$&$0$&$0$\\[0.2cm]
\hline
 &\multicolumn{2}{c}{Family \famc}\\[0.2cm]
&Standard&Non-standard\\[0.1cm]
&$(\xi\ne 0,\vartheta)$
&$(\vartheta,\beta_+,\beta_-)$\\[0.2cm]
\hline
$r$& $\vartheta N\wedge \mmm + \xi A_+\wedge A_-$&
 $(\vartheta N+ \beta_+ A_+ + \beta_- A_-)\wedge \mmm $\\[0.3cm]
$\delta(N)$&$0$&$(\beta_+ A_+ - \beta_- A_-)\wedge \mmm$\\[0.1cm]
$\delta(A_+)$&$-(\vartheta + \xi)A_+\wedge \mmm$
&$- \vartheta  A_+\wedge \mmm $\\[0.1cm]
$\delta(A_-)$&$ (\vartheta - \xi)A_-\wedge \mmm$
&$ \vartheta  A_-\wedge \mmm$\\[0.1cm]
$\delta(\mmm)$&$0$&$0$\\[0.2cm]
 \hline
\end{tabular}}

\bigskip

We recall   that  all  oscillator
bialgebras are coboundary ones  \cite{goslar} and they   were explicitly
obtained in \cite{osc}. In particular:

\noindent
$\bullet$ The family {\fama} corresponds to the type I$_+$ of
\cite{osc} except for the presence of the parameter
$\beta_+$. However this parameter is superfluous: if we
define a new generator as $N'=N+(\beta_+/\alpha_+)\mmm$ we find that the
commutation rules (\ref{nb}) are preserved and  $\beta_+$  appears
explicitly in Table 2.

\noindent
$\bullet$  The
bialgebras of the type I$_-$ of \cite{osc} are completely equivalent to
those of type I$_+$ by means of an automorphism  similar to the one
defined by (\ref{de}) and (\ref{df}) for $\ssl$.

\noindent
$\bullet$  The non-standard family
{\famc} corresponds exactly to the non-standard type II of \cite{osc}
but the standard subfamily does not, that is, the parameters $\beta_+$
and $\beta_-$ do not appear in the contracted bialgebras.
We can introduce them by means of  another
automorphism defined through:
\bea
&&N'=N-\frac{\beta_+}{\vartheta +\xi}A_+ -
\frac{\beta_-}{\vartheta -\xi} A_-\qquad \vartheta +\xi\ne 0
\qquad \vartheta -\xi\ne 0\cr
&&A'_+=A_+ - \frac{\beta_-}{\vartheta -\xi}\mmm \qquad
A'_-=A_- - \frac{\beta_+}{\vartheta +\xi}\mmm \qquad \mmm'=\mmm.
\eea
These new  generators satisfy the commutation rules (\ref{nb})  and
now the standard  family {\famc} can be  identified within the
classification of  \cite{osc}. In particular, the harmonic oscillator
Lie bialgebra corresponding to \cite{GS} is recovered in the case
$\vartheta=0$ and $\xi=-z$. As a byproduct, we have shown that  when
$\vartheta +\xi\ne 0$ and $\vartheta -\xi\ne 0$,   both parameters
$\beta_+$,  $\beta_-$ are irrelevant. 

Therefore, there exist only two isolated
oscillator Lie bialgebras that we do not find by contracting $\ssl$: if
$\vartheta=\xi$ it does not seem possible  to introduce $\beta_-$,
and likewise, if  $\vartheta=- \xi$ to recover $\beta_+$. In the rest
of the cases, the non-relativistic counterparts of $\ssl$
algebraic structures can be easily obtained. In particular, Lie
bialgebra contractions would give rise to (multiparametric) quantum
$h_4$ algebras when applied onto the quantum $\ssl$ deformations
that will be considered in the following Section.


\sect{Multiparametric quantum $\ssl$ algebras}

Now we proceed to obtain some relevant quantum Hopf algebras
corresponding to   the  $\ssl$  bialgebras. We shall write only the
coproducts and the deformed commutation rules as the counit is always
trivial   and the antipode can be easily deduced by means of the Hopf
algebra axioms. We emphasize that coproducts are found by computing a
certain ``exponential" of the Lie bialgebra structure that
characterizes the first order in the deformation. Deformed Casimir
operators, which are essential for the construction of integrable systems,
are also explicitly given.

\subsect{Family \fama\ quantizations}

\subsubsect{Standard subfamily with $\atwo=0$ and $\athree=0$}

If $\atwo$ and $\athree$ vanish, we have that $\aone\ne 0$ and $\afive\ne
0$. Performing the following change of basis
\be
\K'=\K- \frac{\aone}{\afive}\P  ,
\label{za}
\ee
the cocommutator adopts a simpler form
\bea
\delta(\K')=0\qquad \delta(\P)=\afive \K'\wedge \P\qquad
\delta(\T)=\afive \K'\wedge \T\qquad
\delta(\M)= 0 ,
\label{zb}
\eea
while the classical $r$-matrix is formally preserved as
$r=\frac 12(\aone \K'\wedge \P - 2 \afive \P\wedge \T)$.

In this new $\ssl$ basis the commutators (\ref{bc}) and Casimir
(\ref{bd}) turn out to be
\bea
&&[\K',\P]=2\P \qquad [\K',\T]=-2\T -\frac{\aone}{\afive}\K'-
\frac{\aone^2}{\afive^2}\P\cr
&&  [\P,\T]=\K'+\frac{\aone}{\afive}\P \qquad [\M,\cdot\,]=0.
\label{zc}
\eea
\be
{\cal C}=\left(\K'+ \frac{\aone}{\afive}\P\right)^2 + 2 \P \T + 2 \T \P.
\label{zzc}
\ee

The coproduct of the corresponding quantum algebra
$U_{\aone,\afive}(\ssl)$
can be easily deduced from (\ref{zb}) and it reads
\bea
&&\Delta(\K')=1\otimes \K' + \K'\otimes 1\qquad
\Delta(\P) =e^{ \afive \K' /2}\otimes \P + \P\otimes
e^{- \afive \K' /2}\cr
&&\Delta(\M)=1\otimes \M + \M\otimes 1
\qquad
\Delta(\T) =e^{ \afive \K' /2}\otimes \T + \T\otimes
e^{- \afive \K' /2} .
\label{zd}
\eea
We can return to the initial
basis with $\K$ instead of $\K'$; however, in this case it does
not seem worthy since $\K$ and $\P$ do not commute and this fact
would complicate further computations (note also that $\K'$ is
primitive so that we know that $\Delta(e^{x \K'})=e^{x \K'}\otimes
e^{x \K'}$ for any parameter $x$).

Deformed commutation rules compatible with (\ref{zd}) are found to be
\bea
&&[\K',\P]=2\P \qquad [\K',\T]=-2\T -\frac{\aone}{\afive}\,
\frac{\sinh(\afive \K'/2)}{\afive/2}  -
\frac{\aone^2}{\afive^2}\P\qquad [\M,\cdot\,]=0\cr
&&  [\P,\T]=\frac{\sinh \afive \K' }{\afive}
+\frac{\aone}{ \afive} \left(\frac{e^{\afive}-1}{2\afive}\right)
\left(e^{ -\afive \K' /2}\P + \P e^{ \afive \K' /2}\right) ,
\label{ze}
\eea
and the central element that deforms  the Casimir (\ref{zzc}) is
\bea
&&{\cal C}=\frac{2}{\afive\tanh\afive}\left(\cosh(\afive\K') -1\right)+
\frac{\aone}{\afive}\left(\frac{\sinh(\afive \K'/2)}{\afive/2}\P
+\P\frac{\sinh(\afive \K'/2)}{\afive/2}\right)\cr
&&\qquad\quad +\frac{\aone^2}{\afive^2}\P^2+ 2 (\P \T +
\T \P) .
\eea
In order to check these results, the following relations are useful
\bea
&&e^{x\K'}\T e^{-x\K' }=\T e^{- 2 x} +
\frac{\aone}{\afive^2}(e^{-2 x}-1)\sinh(\afive \K'/2)-
\P \frac{\aone^2}{2\afive^2}\sinh 2 x \cr
&&e^{x \K'}\P e^{-x\K'}=\P e^{2 x} .
\eea

We remark that here the $u(1)$ (central) generator $\M$  does not couple
with the $sl(2,\R)$ sector, so
$U_{\aone,\afive}(\ssl)= U_{\aone,\afive}(sl(2,\R))\oplus u(1)$.

It is
also interesting to stress that the (to our knowledge, new)
quantum algebra $U_{\aone,\afive}(sl(2,\R))$ is just a superposition
of the standard and non-standard deformations of
$sl(2,\R)$, since its classical $r$-matrix is the sum of both the
standard and the non-standard one for $sl(2,\R)$.
This fact can be clearly appreciated  by deducing the associated
quantum
$R$-matrix in the fundamental representation.  By following
\cite{BCGST}, we get a $2\times 2$ matrix representation $D$ of
(\ref{ze}) given by
\bea &&
D(\K')=\left(\begin{array}{cc}
1& -\frac{\aone}{\afive} \\ 0&-1
\end{array}\right)\qquad
D(\P)=\left(\begin{array}{cc}
0&\cosh(\frac{\afive}2)\\ 0&0
\end{array}\right)\qquad 
D(\M)=\left(\begin{array}{cc}
1&0 \\ 0&1
\end{array}\right)\cr
&&
D(\T)=\left(\begin{array}{cc}
0&\frac{\aone^2}{4 \afive^2} \left(\frac 2{\afive}\sinh(\frac{\afive}2)-
\cosh(\frac{\afive}2)\right)\\
\frac 2{\afive}\sinh(\frac{\afive}2)&0
\end{array}\right)
\label{repre}
\eea
which, in turn, provides a $4\times 4$ matrix representation of the
coproduct (\ref{zd}). We consider now an arbitrary  $4\times 4$ matrix
and impose it to fulfil both the quantum YBE and the property
\be
{\cal R}\Delta(X){\cal R}^{-1}=\sigma\circ \Delta(X) 
\label{fg}
\ee
for $X\in\{D(\K'),D(\P),D(\T),D(\M)\}$, and where
$\sigma(A\otimes B)=B\otimes A$. Finally we find the solution
\be
{\cal R}=\left(\begin{array}{cccc}
1&h&-qh&h^2\\
0&q&1-q^2&qh\\
0&0&q&-h\\
0&0&0&1
\end{array}\right)
\label{rsn}
\ee
where
\be
q=e^{\afive}\qquad h=\frac{\aone}{2}\left(\frac
{e^{\afive}-1}{\afive}\right) .
\ee

The expression (\ref{rsn}) clearly shows the intertwining between
standard and non-standard properties within the quantum algebra
$U_{\aone,\afive}(\ssl)$. This results in  $\cal R$ being a
quasitriangular solution of quantum YBE and not a triangular one since
${\cal R}_{12}{\cal R}_{21}\ne I$.
In the fundamental representation (\ref{repre}), the standard quantum
$R$-matrix of
$sl(2,\R)$ would be obtained in the limit $\aone\to 0$, and the
non-standard or Jordanian one
\cite{nonsa,nonsb} would be a consequence of taking
$\afive\rightarrow 0$. 
However, we stress that the latter is not a well defined limit at the
Hopf algebra level (see \cite{KLM} for a detailed study of this kind of
problems). Moreover,
$U_{\aone,\afive}(\ssl)$ is just the quantum algebra underlying  the
construction of non-standard quantum 
$R$-matrices out of standard ones proposed in \cite{AKS,ACC} and  its
dual Hopf algebra would give rise to the quantum group $GL_{h,q}(2)$
introduced in
\cite{Kuper}. We finally recall that the classification of
$4\times 4$ constant solutions of the quantum YBE can be found in 
\cite{Hlav}.

\subsubsect{Non-standard subfamily with $\afive=0$}

We restrict to the case with $\afive=0$ so that $\P$ is a primitive
generator.  The coproduct can be easily deduced by applying the
Lyakhovsky--Mudrov  method
\cite{LM} in the same way as in the oscillator $h_4$ case
\cite{osc}.  The cocommutators for the two non-primitive generators can
be written in matrix form as
\be
 \delta\left(\begin{array}{c}
\K' \\ \T
\end{array}\right)=
\left(\begin{array}{cc}
-\aone \P& 0   \\ \frac {\athree}2 \M & -\aone \P
\end{array}\right)\dot\wedge \left(\begin{array}{c}
\K' \\ \T
\end{array}\right)
\label{cp}
\ee
where
\be
\K':= \K-\frac{\athree}{\aone} \M\qquad \aone\ne 0 .
\label{co}
\ee
Hence their  coproduct is given by 
\be
 \Delta\left(\begin{array}{c}
\K' \\ \T
\end{array}\right)=
\left(\begin{array}{ccc}
1 \otimes \K'  \\ 1 \otimes \T
\end{array}\right) +
\sigma\left( \exp\left\{\left(\begin{array}{cc}
 \aone \P & 0  \\ -\frac {\athree}2 \M
  &  \aone \P
\end{array}\right) \right\}\dot\otimes \left(\begin{array}{c}
\K' \\ \T
\end{array}\right)\right)  ,
\label{cq}
\ee
where $\sigma(X\otimes Y):=Y\otimes X$. The exponential
of the Lie bialgebra matrix coming from (\ref{cp}) is the essential
object in the obtention of the deformed coproduct, whose
coassociativity is ensured by construction \cite{LM}.
In terms of the original basis
the coproduct,  commutation rules and  Casimir of the quantum $\ssl$
algebra,   $U_{\aone,\athree}(\ssl)$, are given
by
\bea
&&\Delta(\P)=1\otimes \P + \P\otimes 1\qquad
\Delta(\M)=1\otimes \M + \M\otimes 1\cr
&&\Delta(\K)=1\otimes \K + \K\otimes e^{\aone \P} -
\athree \M\otimes \left(\frac {e^{\aone \P}-1}{\aone}\right) \cr
&&\Delta(\T)=1\otimes \T + \T\otimes e^{\aone \P}
- \frac {\athree}2 \left( \K -
\frac {\athree}{\aone} \M \right)\otimes \M e^{\aone \P} ,
\label{cr}
\eea
\bea
&&[\K, \P] =2 \frac {e^{\aone \P} - 1} {\aone}  \qquad
[\K, \T] = - 2 \T + \frac{\aone}{2}\left(
\K - \frac{\athree}{\aone}\M\right)^2  \cr
&&[\P, \T] = \K +  {\athree} \M  \frac { e^{\aone \P} - 1}{\aone}
 \qquad [\M,\,\cdot \,] = 0,
\label{ct}
\eea
\bea
&&{\cal C}=
\left(\K-\frac{\athree}{\aone} \M\right)
e^{-\aone\P}\left(\K-\frac{\athree}{\aone} \M\right)
+ 2 \frac{\athree}{\aone} \K\M\cr
&&\qquad
+2\frac{1-e^{-\aone\P} }{\aone}\T+
2\T \frac{1-e^{-\aone\P} }{\aone}
+2(e^{-\aone\P}-1).
\eea

It is interesting to note that $U_{\aone,\athree}(\ssl)$
reproduces  the  two-parameter Jordanian deformation of $gl(2)$ obtained
in \cite{Preeti} once we relabel the deformation parameters as
$ \aone = 2 h$ and $\athree = - 2 s$. We also remark that this quantum
deformation has been constructed in \cite{Dobrev} by using a  duality
procedure from the quantum group $GL_{g,h}(2)$ introduced in \cite{Agha};
it can be checked that  the generators
$\{A,B,H,Y\}$ and deformation parameters $\tilde g$, $\tilde h$ defined
by
\bea
&&A=\M\qquad B=\P\cr
&&H=\exp\{-\aone\P/2\} \K + 2 \frac{\athree}{\aone}\M
\sinh(\aone\P/2)\cr
&&Y=\exp\{-\aone\P/2\}\T -
\frac{\athree^2}{4\aone}\exp\{\aone\P/2\}\M^2+
\frac{\aone}{8}\sinh(\aone\P/2)\cr
&&\tilde g= - \aone/2\qquad \tilde h = - \athree/2
\eea
give rise to the quantum $gl(2)$ algebra worked out in \cite{Dobrev}.
On the other hand, we also recover the quantum extended $sl(2,\R)$
algebra introduced in \cite{boson} if we consider the basis
$\{\K',\P,\T,\M\}$ and we set $\aone =2z$ and $\athree=-2z$. The
corresponding universal
$R$-matrix can be also found in \cite{boson,Preeti}.

In the basis here adopted, a coupling of the central
generator
$\M$ with the $sl(2,\R)$ sector arises; however if
we set  $\athree=0$ such coupling disappears and we can rewrite
$U_{\aone}(\ssl)=U_{\aone}(sl(2,\R))\oplus u(1)$
 where  $U_{\aone}(sl(2,\R))$ is the  non-standard or Jordanian
deformation of  $sl(2,\R)$ \cite{nonsa,nonsb,nonsc,nonsd,nonse,nonsf}.

\subsect{Family \famc\ quantizations}

\subsubsect{Standard subfamily and twisted XXZ models}

 The coproduct and commutators  of the two-parametric
 quantum  algebra
$U_{\afive,\asix}(\ssl)$ are given by
\bea
&&\Delta(\M)=1\otimes \M + \M\otimes 1\qquad
 \Delta(\K)=1\otimes \K + \K\otimes 1\cr
&&\Delta(\P)=e^{(\afive \K - \asix \M)/2}\otimes \P + \P\otimes
e^{-(\afive \K - \asix \M)/2} \cr
&&
 \Delta(\T)=e^{(\afive \K + \asix \M)/2}\otimes \T + \T\otimes
e^{-(\afive \K + \asix \M)/2}
\label{er}
\eea
\be
[\K,\P]=2\P \qquad [\K,\T]=-2\T \qquad [\P,\T]=\frac{\sinh
\afive\K}{\afive}
\qquad [\M,\cdot\,]=0.
\label{es}
\ee
The deformed Casimir is
\be
{\cal C}=\cosh \afive \left(\frac{\sinh
(\afive\K/2)}{\afive/2} \right)^2  +2\,\frac{\sinh \afive}{\afive}\,(\P
\T + \T \P).
\label{et}
\ee

This quantum algebra, together with its corresponding universal quantum
$R$-matrix, has been obtained in  \cite{Dobrevb} and   \cite{CJ}, and it
can be related with the so called
$gl_{q,s}(2)$ introduced in \cite{wess} (see also \cite{burdik}) by
defining a set of new generators in the form:
\bea
&&{\widetilde J}_0 =\frac 12 \K \qquad
{\widetilde J}_+ =\sqrt{\frac{\afive}{\sinh \afive}}
\exp\left\{ \frac{\asix} {2\afive}\M \right\}\P \cr
&&
{\widetilde Z} = \frac{\asix} {2\afive}\M \qquad
{\widetilde J}_- =\sqrt{\frac{\afive}{\sinh \afive}}
\exp\left\{ -\frac{\asix} {2\afive}\M \right\}\T
 \qquad
\label{is}
\eea
and the parameters $q$ and $s$ as
\be
q=e^\eta\qquad \eta=-\afive\qquad s^{ {\tilde Z}}=
\exp\left\{ \frac{\asix} {2\afive}\M \right\} .
\label{isp}
\ee

The algebra
$U_{\afive,\asix}(\ssl)$ is just the quantum algebra underlying the
XXZ Heisenberg Hamiltonian with twisted boundary conditions
\cite{MRplb}.  This
deformation can be thought as a Reshetikhin twist of the usual
standard deformation. This superposition of the standard
quantization and a twist is easily reflected at the Lie bialgebra
level by the associated classical $r$-matrix $r=- \frac 12  \asix
\K\wedge \M - \afive \P\wedge\T$ (see Table 1): within it, the second
term generates the standard deformation and the exponential of the
first one gives us the Reshetikhin twist. Compatibility between both
quantizations is ensured by the fact that $r$ fulfils the
modified classical YBE and the method here used shows that the full
simultaneous quantization of the two-parameter Lie bialgebra is
possible.

On the other hand, if we set $\asix=0$  we find that $\M$ does not
couple with the deformation of the $sl(2,\R)$ sector, and
$U_{\afive}(\ssl)=U_{\afive}(sl(2,\R))\oplus u(1)$
 where  $U_{\afive}(sl(2,\R))$ is the well-known standard deformation of
$sl(2,\R)$ \cite{standb,standa}.

\subsubsect{Non-standard subfamily and twisted XXX models}

 This bialgebra has one primitive generator    $\M$;
the cocommutator  for the  remaining generators can be written as
\be
 \delta\left(\begin{array}{c}
\K \\ \P\\ \T
\end{array}\right)=
\left(\begin{array}{ccc}
0 & -\athree \M & -\afour \M \\ \frac{ \afour}2 \M& -\asix \M
& 0\\
\frac {\athree}2 \M &0& \asix \M
\end{array}\right)\dot\wedge \left(\begin{array}{c}
\K \\ \P\\ \T
\end{array}\right) ,
\label{fa}
\ee
so that their coproduct is given by:
\be
 \Delta\left(\begin{array}{c}
\K \\ \P\\ \T
\end{array}\right)=
\left(\begin{array}{ccc}
1 \otimes \K  \\ 1 \otimes \P \\ 1 \otimes \T
\end{array}\right) +
\sigma \left( \exp\left\{\left(\begin{array}{ccc}
 0 & \athree \M & \afour \M \\ -\frac{ \afour}2 \M & \asix \M
& 0\\
 -\frac {\athree}2 \M&0&  -\asix \M
\end{array}\right)\right\}\dot\otimes \left(\begin{array}{c}
\K \\ \P\\ \T
\end{array}\right)\right) .
\label{fb}
\ee
 If we denote
the exponential of the Lie bialgebra matrix by $E$,  the coproduct
can be expressed in terms of the $E_{ij}$ entries as follows:
\bea
&&\Delta(\M)=1\otimes \M + \M \otimes 1 \cr
&&\Delta(\K)=1\otimes \K + \K \otimes E_{11}(\M) +
\P \otimes E_{12}(\M) + \T \otimes E_{13}(\M) \cr
&&\Delta(\P)=1\otimes \P + \P \otimes E_{22}(\M) +
\K \otimes E_{21}(\M) + \T \otimes E_{23}(\M) \cr
&&\Delta(\T)=1\otimes \T + \T \otimes E_{33}(\M) +
\K \otimes E_{31}(\M) + \P \otimes E_{32}(\M) .
\label{ffe}
\eea
The explicit form of the functions $E_{ij}$ is quite complicated, which
in turn makes it difficult to find the associated deformed  commutation
relations.  Therefore in the sequel we study a specific case by setting
$\afour=0$. The coproduct of the quantum algebra
$U_{\athree,\asix}(\ssl)$ is
\bea
&&\Delta(\M)=1\otimes \M + \M \otimes 1\qquad
\Delta(\P)=1\otimes \P + \P \otimes e^{\asix \M}\cr
&&\Delta(\K)=1\otimes \K + \K \otimes 1 + \athree \P\otimes \left(\frac
{e^{\asix \M}-1}{\asix}\right)\cr
&&\Delta(\T)=1\otimes \T + \T \otimes e^{-\asix \M} + \athree \K\otimes
\left(\frac {e^{-\asix \M}-1}{2\asix}\right) \cr
&&\qquad \qquad+ \athree^2 \P\otimes
\left(\frac {1 - \cosh {\asix \M}}{2\asix^2}\right)
\label{fd}
\eea
and the associated commutation rules are the non-deformed ones
(\ref{bc}). The role of $\M$ is
essential in this deformation, and no uncoupled structure can be
recovered unless all deformation parameters vanish. On the other hand,
the element
\be
{\cal R}= \exp\{r\}=
\exp\{\M\otimes (\asix\K -\athree \P)/2\}
\exp\{-(\asix\K -\athree \P)\otimes \M /2\}
\label{fe}
\ee
is a solution of the quantum YBE (as $\M$ is a central generator) and it
also fulfils the relation (\ref{fg}). The proof of this property is
sketched in the Appendix. In the fundamental representation
(\ref{frep}), the $R$-matrix (\ref{fe}) reads
\be
{\cal R}=\left(\begin{array}{cccc}
1&-e^{-\asix}\,p&p&-e^{-\asix}\,p^2\\
0&e^{-\asix}&0&e^{-\asix}\,p\\
0&0&e^{\asix}&-p\\
0&0&0&1
\end{array}\right)\qquad
p=\frac{\athree}{2}\left(\frac{e^{\asix}-1}{ \asix}\right).
\label{rxxt}
\ee

From the point of view of spin systems, a direct connection between the
one-parameter deformation with $\athree=\afour=0$ and the
twisted XXX chain can be established. This particular quantization can be
obtained as the limit
$\afive\to 0$ of the (standard) quantum algebra
$U_{\afive,\asix}(\ssl)$, and $\afive$ is known to be related to the
anisotropy of the XXZ model. Under such a limit, twisted boundary
conditions coming from $\asix\neq 0$ are preserved, and a twisted XXX
model is expected to arise. This symmetry property can be explicitly
checked by following the approach presented in \cite{ACF} (and used there
in order to obtain deformed t-J models). We consider the fundamental
representation $D$ of $U_{\asix}(\ssl)$ (\ref{frep})
in terms of Pauli spin matrices:
$D(\K)=\sigma_3$, $D(\P)=\sigma_+$, $D(\T)=\sigma_-$ and $D(\M)$
will be again the two-dimensional identity matrix. If we compute
(with $\athree=0$) the deformed coproduct (\ref{fd})  of the
Casimir (\ref{bd}), we obtain
\be
(D\otimes D)(\Delta_{\asix}({\cal{C}}))= 6 + 2\, ( 
\sigma_3\otimes \sigma_3 + 2\,e^{-\asix}\,\sigma_-\otimes \sigma_+ +
2\,e^{\asix}\sigma_+\otimes \sigma_-).
\label{txxx}
\ee 
This means that the twisted XXX Heisenberg Hamiltonian can be written (up
to global constants) as the sum of elementary two-site Hamiltonians
given by the coproducts (\ref{txxx}) 
\bea
&& H_{\asix}=\sum_{i=1}^{N}(D_i\otimes
D_{i+1})(\Delta_{\asix}^{i,i+1}({\cal{C}}))\cr
&& \qquad=6\,N   + 2\, \sum_{i=1}^{N}{ ( 
\sigma_3^i \sigma_3^{i+1} + 2\,e^{-\asix}\,\sigma_-^i \sigma_+^{i+1} +
2\,e^{\asix}\sigma_+^i \sigma_-^{i+1}) }.
\label{txxxh}
\eea
This expression explicitly reflects the $U_{\asix}(\ssl)$ 
quantum algebra invariance of this model since, by construction, the
Hamiltonian (\ref{txxxh}) commutes with the $(N+1)$-th coproduct of the
generators of $U_{\asix}(\ssl)$. In the same way, further contributions
could be obtained by considering other quantum deformations belonging to
this family. In particular, if we take the two-parametric coproduct
(\ref{fd}) and repeat the same construction we are lead to the following
spin Hamiltonian
\bea
&& \!\!\!\!H_{\athree,\asix}=\sum_{i=1}^{N}(D_i\otimes
D_{i+1})(\Delta_{\athree,\asix}^{i,i+1}({\cal{C}})) \cr
&& \qquad=6\,N  + 2\, \sum_{i=1}^{N} {\big( 
\sigma_3^i \sigma_3^{i+1} + 2\,e^{-\asix}\,\sigma_-^i \sigma_+^{i+1} +
2\,e^{\asix}\sigma_+^i \sigma_-^{i+1}\big) }\cr
&& \quad\qquad + 2\athree 
\sum_{i=1}^{N} \left\{ \left(\frac
{e^{\asix}-1}{\asix}\right)\,\sigma_+^i\sigma_3^{i+1}
+ \left(\frac
{e^{-\asix}-1}{\asix}\right)\,\sigma_3^i\sigma_+^{i+1}
\right\} \cr
&& \qquad\qquad + 2 \athree^2 \sum_{i=1}^{N} \left(\frac
{1-\cosh{\asix}}{\asix^2}\right)\,\sigma_+^i\sigma_+^{i+1}.
\label{txdef}
\eea
Therefore, we have obtained a (quadratic in $\athree$) deformation of the
twisted XXX chain, which is invariant under $U_{\athree,\asix}(\ssl)$ and
whose associated quantum
$R$-matrix would be (\ref{rxxt}). Likewise, the introduction of the
full quantization containing $\afour$ would provide a further deformation
of the Hamiltonian (\ref{txdef}).

\sect{Concluding remarks}

We have presented a constructive overview of
multiparameter quantum $gl(2)$ deformations based on the
classification and further quantization of $gl(2)$ Lie bialgebra
structures. The quantization procedure (based on a the construction
of the ``exponential" of the first order of the deformation) turns
out to be extremely efficient in order to construct explicitly
multiparametric quantum
$gl(2)$ algebras. By following this method, a family of
new multiparametric quantizations generalizing the symmetries  of twisted
XXX models is introduced and the quantum algebra counterpart of the 
superposition of standard and non-standard deformations $GL_{h,q}(2)$
\cite{Kuper} is obtained. 
Throughout the paper, Lie bialgebra analysis is
shown to provide  essential algebraic information
characterizing the quantum algebras and their associated
models. For instance, a Lie bialgebra contraction method gives
a straightforward way to implement the non-relativistic limit
of the quantum
$gl(2)$ algebras, different coupling possibilities between the central
generator and the $sl(2,\R)$ substructure are easily extracted
from the cocommutator
$\delta$, and Reshetikhin twists giving rise to twisted XXZ
models can be identifyied (and explicitly constructed) with the
help of the classical
$r$-matrices generating the Lie bialgebras. In general, we can conclude
that the existence of a central generator strongly increases   the number
of different quantizations even when this central extension is
cohomologically trivial at the non-deformed level (compare the
classification here presented with the one corresponding  to
$sl(2,\R)$), and the explicit construction of these
quantizations provide an algebraic background for the
systematic obtention of new integrable systems.


\bigskip
\bigskip

\noindent
{\Large{{\bf Acknowledgments}}}

\bigskip

 A.B. and
F.J.H. have been partially supported by DGICYT (Project  PB94--1115)
from the Ministerio de Educaci\'on y Cultura  de Espa\~na and by Junta
de Castilla y Le\'on (Projects CO1/396 and CO2/297). P.P. has been
supported by a fellowship from AECI, Spain and is also grateful to Prof.
Enrique \'Alvarez and Prof. C\'esar G\'omez for providing hospitality at
Dept. de F\1sica Te\'orica, Universidad Aut\'onoma de Madrid.

\bigskip
\bigskip


\noindent
{\Large{{\bf Appendix}}}

\appendix

\setcounter{equation}{0}

\renewcommand{\theequation}{A.\arabic{equation}}

\bigskip

\noindent
We prove that the element (\ref{fe}) satisfies (\ref{fg}) for the
generator $\T$.  We write the universal $R$-matrix as
${\cal R}=\exp\{\M\otimes A\}
\exp\{-A\otimes \M \}$ where $A=\frac 12 (\asix\K -\athree \P)$ and we
take into account the formula
\be
e^{f}\,\Delta(X)\,e^{-f}=\Delta(X) +\sum_{n=1}^\infty \frac
1{n!}\,[f,\dots [f,\Delta(X)]^{n)}\dots].
\label{apenda}
\ee
We set  $f\equiv -A\otimes\M$ and we consider the
coproduct of $\T$ (\ref{fd}), thus we obtain
\bea
&&[f,\dots [f,\Delta(\T)]^{n)}\dots]=
 (\T  +\frac{\athree}{2\asix}\K)\otimes  (\asix I)^n e^{-\asix I} \cr
&&\qquad +\frac{\athree^2}{2\asix^2}\P\otimes
 (\asix I)^n \sinh \asix I\quad \mbox{for  $n$ odd and $n\ge 1$} ;\\
&&[f,\dots [f,\Delta(\T)]^{n)}\dots]=
 (\T  +\frac{\athree}{2\asix}\K) \otimes  (\asix I)^n e^{-\asix I} \cr
&&\qquad -\frac{\athree^2}{2\asix^2}\P\otimes
 (\asix I)^n \cosh \asix I\quad \mbox{for  $n$ even and $n\ge 2$} .
\nonumber
\eea
Therefore,
\bea
&&
e^{f}\Delta(\T)e^{-f}=\Delta(\T)+ (\T +
\frac{\athree}{2\asix}\K )
\otimes  e^{-\asix I}\sum_{n=1}^\infty \frac{(\asix I)^n}{n!}\cr
&&\qquad + \frac{\athree^2}{2\asix^2}\P\otimes \sinh\asix \M
\sum_{k=0}^\infty \frac{(\asix I)^{2k+1}}{(2k+1)!}
- \frac{\athree^2}{2\asix^2}\P\otimes \cosh\asix \M
\sum_{k=1}^\infty \frac{(\asix I)^{2k}}{(2k)!}\cr
&&=\Delta(\T)+ (\T +
\frac{\athree}{2\asix}\K )
\otimes (1 -e^{-\asix I})
+ \frac{\athree^2}{2\asix^2}\P\otimes ( \cosh \asix \M -1)\cr
&&=1\otimes \T + \T\otimes 1\equiv \Delta_0(\T).
\label{apendb}
\eea
Now we take $f\equiv \M\otimes A$ and we find that
\bea
&&[f,\dots [f,\Delta_0(\T)]^{n)}\dots]=
 - (\asix I)^n  \otimes (\T +\frac{\athree}{2\asix}\K)
\quad \mbox{for  $n$ odd and $n\ge 1$} \\
&&[f,\dots [f,\Delta_0(\T)]^{n)}\dots]=
  (\asix I)^n  \otimes (\T +\frac{\athree}{2\asix}\K -
\frac{\athree^2}{2\asix^2}\P)  \quad \mbox{for  $n$ even and $n\ge 2$} .
\nonumber
\eea
And finally the proof follows from
\bea
&&e^{f}\Delta_0(\T)e^{-f}= 1\otimes \T + \T\otimes 1 - \sinh\asix \M
\otimes (\T +\frac{\athree}{2\asix}\K) \cr
&&\qquad  +(\cosh\asix I -1)\otimes (\T
+\frac{\athree}{2\asix}\K  - \frac{\athree^2}{2\asix^2}\P) = \sigma\circ
\Delta(\T) .
\eea
Likewise, it can be checked that (\ref{fg}) is fulfilled for the remaining
generators.

\end{document}